\documentclass[3p,times]{elsarticle}

\usepackage{amsmath}

\usepackage{amssymb}
\usepackage[figuresright]{rotating}

\usepackage[english]{babel}
\begin{document}
\begin{frontmatter}
\large 

\title{Hermite Calculus}

\author[Enea]{G. Dattoli}
\ead{giuseppe.dattoli@enea.it}

\author[LaS]{B. Germano}
\ead{bruna.germano@sbai.uniroma1.it}

\author[Enea,Unict]{S. Licciardi\corref{cor} }
\ead{silvia.licciardi@dmi.unict.it}

\author[LaS]{M. R. Martinelli}
\ead{martinelli@dmmm.uniroma1.it}

\address[Enea]{ENEA - Frascati Research Center, Via Enrico Fermi 45, 00044, Frascati, Rome, Italy}
\cortext[cor]{Corresponding author}
\address[LaS]{University of Rome, La Sapienza, Department of Methods and Mathematic Models for Applied Sciences, Via A. Scarpa, 14, 00161 Rome, Italy}
\address[Unict]{University of Catania, Department of Mathematics, Via Santa Sofia, 64, 95125 Catania, Italy}
\address[LaS]{University of Rome, La Sapienza, Department of Methods and Mathematic Models for Applied Sciences, Via A. Scarpa, 14, 00161 Rome, Italy}

\begin{abstract}
We develop a new method of umbral nature to treat blocks of Hermite and of Hermite like polynomials as independent algebraic quantities. The Calculus we propose allows the formulation of a number of "practical rules" allowing significant simplifications in computational problems. 
\end{abstract}


\end{frontmatter}


\large 

In this letter we deal with a protocol, which will be referred as Hermite calculus, useful to treat computations involving Hermite polynomials and their generalizations as well.\\

 To  give a flavour of the techniques we will employ we consider the integral

\begin{equation}\label{int1}
I(\alpha,\beta,\gamma)=\int_{-\infty}^{\infty}e^{-(\alpha+\beta)x^{2}-\gamma x}dx
\end{equation}
which can be evaluated with ordinary means, thus getting

\begin{equation}
I(\alpha,\beta,\gamma)=\sqrt{\dfrac{\pi}{\alpha + \beta}}e^{\frac{\gamma^{2}}{4(\alpha + \beta)}}
\end{equation}

We will test the formalism, we are going to decribe, by using such a benchmark and restyle eq. \eqref{int1} as it follows

\begin{equation}\label{htilde}
I(\alpha,\beta,\gamma)=\int_{-\infty}^{\infty}e^{-\alpha x^{2}-\hat{h}_{(\gamma.-\beta)}x}dx
\end{equation}
where we have introduced the notation

\begin{equation}\label{serieh}
e^{-\hat{h}_{(\gamma.-\beta)}x}=\sum_{r=0}^{\infty}\dfrac{(-x)^{r}}{r!}\hat{h}_{(\gamma.-\beta)}^{r}=\sum_{r=0}^{\infty}\dfrac{(-x)^{r}}{r!}H_{r}(\gamma,-\beta)
\end{equation}
based on the use of the umbral identity  (\cite{Eric}, \cite{G.Dattoli})\footnote{We have used the umbral notation in a rather inaccurate way, without specifying on which space the operators are acting.
For further comments and an appropriate discussion of the formal content, see the second of refs. \cite{G.Dattoli}.}

\begin{equation}\begin{split}\label{herm}
& \hat{h}_{(\gamma.-\beta)}^{r}=H_{r}(\gamma,-\beta),\\
& H_{r}(x,y)=r!\sum_{s=0}^{\left[ \frac{r}{2}\right] }\dfrac{x^{r-2s}y^{s}}{(r-2s)!s!}
\end{split}\end{equation}

In the integral in eq. \eqref{serieh} we have treated the term which can be expended in terms of Hermite polynomials as a single block and we have enucleated the variable x raised to the first power.\\

We will now follow a prescription according to which the operator $\hat{h}$  is treated as an ordinary algebraic quantity. \\
According to the ordinary rules for the Gaussian integrals we can write \cite{G.Dattoli} \cite{H.M.Srivastava}

\begin{equation}
I(\alpha,\beta,\gamma)=\sqrt{\dfrac{\pi}{\alpha}}e^{\frac{\hat{h}_{(\gamma.-\beta)}^{2}}{4\alpha}}=\sqrt{\dfrac{\pi}{\alpha}}\sum_{r=0}^{\infty}\dfrac{1}{r!}\left(\dfrac{\hat{h}_{(\gamma.-\beta)}^{2}}{4\alpha} \right)^{r} 
\end{equation}
which provides us with the correct result for the problem we are studying.  The application of the previous prescription yields, indeed, if $\left|  \dfrac{\beta}{\alpha}\right|  < 1$

\begin{equation}\begin{split}\label{funcgen}
& \sqrt{\dfrac{\pi}{\alpha}}\sum_{r=0}^{\infty}\dfrac{1}{r!}\left(\dfrac{\hat{h}_{(\gamma.-\beta)}^{2}}{4\alpha} \right)^{r}=\\
& =\sqrt{\dfrac{\pi}{\alpha}}\sum_{r=0}^{\infty}\dfrac{1}{r!}H_{2r}\left( \dfrac{\gamma}{2\sqrt{\alpha}},-\dfrac{\beta}{4\alpha}\right)=\sqrt{\dfrac{\pi}{\alpha+\beta}}e^{\frac{\gamma^{2}}{4(\alpha+\beta)}} 
\end{split}\end{equation}
which is obtained after using the identity \cite{G.Dattoli}

\begin{equation}
\sum_{r=0}^{\infty}\dfrac{t^{r}}{r!}H_{2r}(x,y)=\dfrac{1}{\sqrt{1-4yt}}e^{\frac{x^{2}t}{1-4yt}}
\end{equation}
The result in eq. \eqref{funcgen}, evidently correct,  yields some confidence on the reliability of the formalism, which is based on a rather wild use of the umbral formalism.
\\

Even though not explicitly stated,  the umbral operator defined in eq. \eqref{herm} satisfies the identity\footnote{The subscript $\left( \gamma,-\beta\right) $ has been omitted because the identity holds for $\hat{h}$ operators with the same basis, hereafter it will be included whenever necessary. }

\begin{equation}
\hat{h}^{m}\hat{h}^{r}=\hat{h}^{m+r}
\end{equation}
where the powers $m$ and $r$ are not necessarily real integers. It is also fairly natural to set

\begin{equation}
\partial_{\hat{h}}\hat{h}^{r}=r\;\hat{h}^{r-1}=rH_{r-1}(\gamma,-\beta)
\end{equation}
Since the following recurrence holds

\begin{equation}
\partial_{\gamma}H_{r}(\gamma,-\beta)=rH_{r-1}(\gamma,-\beta)
\end{equation}
the "derivative" operator can therefore be identified with

\begin{equation}
\partial_{\hat{h}}\rightarrow \partial_{\gamma}
\end{equation}
Furthermore since

\begin{equation}
\hat{h}\hat{h}^{r}=\hat{h}^{r+1}=H_{r+1}(\gamma,-\beta)
\end{equation}
and, on account of the recurrence,

\begin{equation}
H_{r+1}(\gamma,-\beta)=\gamma H_{r}(\gamma,-\beta)-2\; \beta \;r\; H_{r-1}(\gamma,-\beta)
\end{equation}
we can also conclude that $\hat{h}$ itself can be identified with the differential operator

\begin{equation}\label{opPDE}
\hat{h}=\gamma -2\; \beta \;\partial_{\gamma}
\end{equation}

It is also worth noting that

\begin{equation}\begin{split}
& \partial_{x}^{r}e^{-\hat{h}x}=(-1)^{r}\hat{h}^{r}e^{-\hat{h}x}=(-1)^{r}\sum_{n=0}^{\infty}\dfrac{(-x)^{n}}{n!}\hat{h}^{n+r}=\\
& =(-1)^{r}\sum_{n=0}^{\infty}\dfrac{(-x)^{n}}{n!}H_{n+r}(\gamma,-\beta)
\end{split}\end{equation}
and, according to the identity

\begin{equation}
\sum_{n=0}^{\infty}\dfrac{t^{n}}{n!}H_{n+l}(x,y)=H_{l}(x+2yt,y)e^{xt+yt^{2}}
\end{equation}
we can establish the “rule”

\begin{equation}\label{rule}
\hat{h}^{r} e^{-\hat{h}x}=H_{r}(\gamma +2\;\beta \;x,-\beta)e^{-(\gamma x +\beta x^{2})}
\end{equation}

We can now make a step further by defining the integral

\begin{equation}\begin{split}
& I(\gamma,\beta)=\int_{-\infty}^{\infty}e^{-\hat{h}x^{2}}dx\\
&e^{-\hat{h}x^{2}}= e^{-(\gamma x^{2}+\beta x^{4})}
\end{split}\end{equation}
which, after applying the prescription that $\hat{h}$ can be treated as an ordinary algebraic quantity, writes

\begin{equation}\label{Ih}
I(\gamma,\beta)=\sqrt{\pi}\hat{h}^{-\frac{1}{2}}
\end{equation}
which makes sense only if we can provide a meaning for $\hat{h}^{-\frac{1}{2}}$ , the most natural conclusion is that they can be understood as factional order Hermite, which for our purposes can be defined as it follows \cite{R.Hermann}

\begin{equation}\label{Hermfrazcos}
H_{\nu}(x,-y)=y^{\frac{\nu}{2}}\dfrac{e^{ \frac{x^{2}}{4y}}}{\sqrt{\pi}}\int_{0}^{\infty}e^{-\frac{t^{2}}{4}}t^{\nu}\cos \left( \dfrac{x}{2\sqrt{y}}t-\dfrac{\pi}{2}\nu\right)dt 
\end{equation}
or as

\begin{equation}\label{Hermfraz}
H_{\nu}(x,-y)=\Gamma(\nu+1)\sum_{r=0}^{\infty}\dfrac{x^{\nu-2r}(-y)^{r}}{\Gamma(\nu+1-2r)r!}, \qquad x>>y
\end{equation}
which has however a limited range of convergence.\\
 The correctness of eq. \eqref{Ih} can be readily proved by a numerical check, involving either the definitions \eqref{Hermfrazcos} and \eqref{Hermfraz}.\\
We will comment, later in this paper, on the extension of the Hermite polynomials to non-integer index.\\
 
Let us now consider the following repeated derivatives

\begin{equation}\begin{split}
& \partial_{x}^{\; n}e^{-\hat{h}x^{2}}=(-1)^{n}H_{n}(2\;\hat{h}\;x,-\hat{h})e^{-\hat{h}x^{2}}=\\
& = (-1)^{n}n!\sum_{r=0}^{\left[ \frac{n}{2}\right] }\dfrac{(-1)^{r}(2x)^{n-2r}}{(n-2r)!r!}\left( \hat{h}^{n-r}e^{-\hat{h}x^{2}}\right) 
\end{split}\end{equation}\\

Thus, getting, on account of eq. \eqref{rule},

\begin{equation}
\partial_{x}^{\;n}e^{-\hat{h}x^{2}}=(-1)^{n}n!\sum_{r=0}^{\left[ \frac{n}{2}\right] }\dfrac{(-1)^{r}(2x)^{n-2r}}{(n-2r)!r!}H_{n-r}(\gamma +2\;\beta x^{2},-\beta)e^{-(\gamma x^{2} +\beta x^{4})}
\end{equation}
in accordance with

\begin{equation}
\partial_{x}^{\;n}e^{-(\gamma x^{2}+\beta x^{4})}=H_{n}^{(4)}(-2\;\gamma x-4\;\beta x^{3},-\gamma-6\;\beta x^{2},-4\;\beta x,-\beta)e^{-(\gamma x^{2}+\beta x^{4})}
\end{equation}

In ref. \cite{J.Bohacik} the following integral

\begin{equation}\begin{split}\label{intJ}
& J(a,b,c) = \int_{-\infty}^{\infty} e^{-(ax^{4}+bx^{2}+cx)}dx,\\
& Re(a)>0
\end{split}\end{equation} 
has been considered, within the framework of problems regarding the non-perturbative treatment of the anharmonic oscillator. A possible perturbative treatment is that of setting

\begin{equation}\begin{split}
& J(a,b,c) = \sum_{n=0}^{\infty}\dfrac{(-1)^{n}}{n!}g_{n}(a)\left[ H_{n}(c,-b)+H_{n}(-c,-b)\right] ,\\
& g_{n}(a) = \int_{0}^{\infty}x^{n}e^{-ax^{4}}dx=\frac{1}{4}a^{-\frac{n+1}{4}}\Gamma\left( \dfrac{n+1}{4}\right) 
\end{split}\end{equation}
which, as noted in \cite{J.Bohacik}, is an expansion with zero radius of convergence in spite of the fact that $J(a,b,c)$ is an entire function for any real or complex value of $b, c$.\\

The use of our point of view allows to write

\begin{equation}\label{Jabc}
J(a,b,c)=\int_{-\infty}^{\infty}e^{-\hat{h}_{(b,-a)}x^{2}-cx}dx=\sqrt{\dfrac{\pi}{\hat{h}}}e^{\frac{c^{2}}{4\hat{h}}}=\sqrt{\pi}\sum_{s=0}^{\infty}\dfrac{1}{s!}\left( \dfrac{c}{2}\right)^{2s}\hat{h}^{-\left( s+\frac{1}{2}\right) } 
\end{equation}
We have omitted the subscript $(b,-a)$ in the r.h.s. of eq. \eqref{Jabc} to avoid a cumbersome notation. The meaning of the operator $\hat{h}$ raised to a negative exponent is easily understood as

\begin{equation}
\hat{h}^{-\left( s+\frac{1}{2}\right) }=H_{-\left( s+\frac{1}{2}\right)}(b,-a)
\end{equation}
where the negative index Hermite polynomials are expressed in terms of the parabolic cylinder functions $D_{n}$ according to the identity \cite{Abramovitz}

\begin{equation}\label{pcf}
 H_{-n}(x,-y)= (2y)^{-\frac{n}{2}}e^{\frac{x^{2}}{8y}} D_{-n}\left( \dfrac{x}{\sqrt{2y}}\right)
\end{equation}

The use of eq. \eqref{pcf} in eq. \eqref{Jabc} finally yields the same series expansion obtained in ref. \cite{J.Bohacik}

\begin{equation}\label{Jserie}
J(a,b,c)=\sqrt{\pi}\sum_{s=0}^{\infty}\dfrac{1}{s!}\left( \dfrac{c}{2}\right)^{2s}\left(2a \right)^{-\frac{1}{2}\left( s+\frac{1}{2}\right) }e^{\frac{b^{2}}{8a}}D_{-\left( s+\frac{1}{2}\right) }\left( \dfrac{b}{\sqrt{2a}}\right)   
\end{equation}
which is convergent for any value of $b, c$ and $a>0$.\\

Regarding the use of non-integer  Hermite polynomials it is evident that the definition adopted in eq. \eqref{Hermfrazcos} can be replaced by the use of the parabolic cylinder fuction, it is therefore worth noting that the use of the properties of the $D$ functions allows the following alternative form for eq. \eqref{Ih} (see ref. \cite{Weisstein})

\begin{equation}\begin{split}
 I(\gamma,\beta) & = \sqrt{\pi} \left(2\;\beta \right)^{-\frac{1}{4}}e^{\frac{\gamma^{2}}{8\beta}}D_{-\frac{1}{2}}\left(\dfrac{\gamma}{\sqrt{2\beta}} \right)= \\
& =\sqrt{\dfrac{\gamma}{2\sqrt{2\beta}}} \left(2\;\beta \right)^{-\frac{1}{4}}e^{\frac{\gamma^{2}}{8\beta}} K_{\frac{1}{4}}\left(\dfrac{\gamma^{2}}{8\beta} \right) ;\\
& D_{-\frac{1}{2}}(z)=\sqrt{\dfrac{z}{2\pi}}K_{\frac{1}{4}}\left( \dfrac{1}{4}z^{2}\right) 
\end{split}\end{equation}
where $K_{\nu}(z)$ is a modified Bessel function of the second kind.\\

A further example of application of the method developed so far is provided by

\begin{equation}\begin{split}
& \int_{0}^{\infty}e^{-(\beta x^{2n}+\gamma x^{n})}dx=\int_{0}^{\infty}e^{-\hat{h}_{(\gamma,-\beta)}x^{n}}dx=\dfrac{1}{n}\Gamma\left(\dfrac{1}{n} \right) \hat{h}_{(\gamma,-\beta)}^{-\frac{1}{n}}=\\
& =\dfrac{1}{n}\Gamma\left(\dfrac{1}{n} \right)\left(2\beta \right)^{-\frac{1}{2n}}e^{\frac{\gamma^{2}}{8\beta}}D_{-\frac{1}{n}}\left( \dfrac{\gamma}{\sqrt{2\beta}}\right) 
\end{split}\end{equation}

In a forthcoming more detailed note we will extend the method to cases involving higher order Hermite polynomials. Just to provide an idea of the extension of the technique, we note that the use of this family of polynomials allows to cast the integral in eq. \eqref{Jabc} in the form

\begin{equation}\label{corr}
J(a,b,c)=\int_{-\infty}^{\infty}e^{\;{}_{4}\hat{h}_{(-c,-a)}x-bx^{2}}dx=\sqrt{\dfrac{\pi}{b}}e^{\frac{\left(\;{}_{4}\hat{h}_{(-c,-a)} \right)^{2} }{4b}}
\end{equation}
where

\begin{equation}\label{4h}
\left(\;{}_{4}\hat{h}_{(-c,-a)} \right)^{n}=H_{n}^{(4)}(-c,-a)=(-1)^{n}n!\sum_{r=0}^{\left[ \frac{n}{4}\right] }\dfrac{c^{n-4r}(-a)^{r}}{(n-4r)!r!}
\end{equation}
with $H_{n}^{(4)}(c,-a)$ being a fourth order Hermite Kamp\'e de F\'eri\'et \cite{P.Appél} polynomial.\\
 The series expansion of the right hand side of eq. \eqref{corr} in terms of fourth order Hermite  converges in a much more limited range than the series \eqref{Jserie} and has been proposed to emphasize the possibilities of the method we have proposed so far.
\\

According to our formalism the Pearcey integral, widely studied in optics, within the framework of diffraction problems \cite{Lopez}, is easily reduced to a particular case of eq. \eqref{intJ}, namely

\begin{equation}
J(1,x,-iy)=\int_{-\infty}^{\infty}e^{-(t^{4}+xt^{2})+iyt}dt=\sqrt{\dfrac{\pi}{\hat{h}_{(x,-1)}}}e^{-\frac{y^{2}}{4\hat{h}_{(x,-1)}}}
\end{equation}

and can be expressed in terms of parabolic cylinder functions, as indicated before. It is perhaps worth stressing that, in the literature a converging  series for the Pearcey integral is given in the form \cite{Berry}

\begin{equation}\begin{split}
& J(1,x,-iy)=\int_{-\infty}^{\infty}e^{-t^{4}-\hat{h}_{(iy,-x)}t}dt=\int_{0}^{\infty}e^{-t^{4}}\left( e^{\hat{h}_{(iy,-x)}t}+e^{\hat{h}_{(-iy,-x)}t}\right) dt=\\
& =2\sum_{n=0}^{\infty}(-1)^{n}g_{2n}(1)a_{2n}(x,y)
\end{split}\end{equation}
with

\begin{equation}\begin{split}
& a_{0}(x,y)=1,\\
& a_{1}(x,y)=y,\\
& a_{n}(x,y)=\dfrac{1}{n}\left(y\;a_{n-1}(x,y)+2\;x\;a_{n-2}(x,y) \right) 
\end{split}\end{equation}
which is reconciled with our previous result, in terms of two variable Hermite polynomials, provided that one recognizes

\begin{equation}
J(1,x,-iy)=\sum_{n=0}^{\infty}\dfrac{(-1)^{n}}{n!}g_{n}(1)\left[ H_{n}(-iy,-x)+H_{n}(iy,-x)\right] 
\end{equation}

In this letter we have provided some hint on the use of the Hermite calculus to study integral forms with specific application in different field of research. In a forthcoming more detailed note we will show how the method can be extended to a systematic investigation of the Voigt functions and  to the relevant generalizations \cite{Pathan}.\\

Before closing the paper we want to underline that the possibilities for the applicability of the integration method discussed in this letter arise if, inside the integrand, an exponential generating function is recognized. \\

To clarify this point we note that the integral

\begin{equation}\begin{split}\label{A}
& f(a,b,c)=\int_{-\infty}^{\infty}e^{-ax^{2}+\sqrt{x^{2}+bx+c}}dx,\\
& b^{2}-4c<0\\
& a>1
\end{split}\end{equation}
can be written as

\begin{equation}\label{B}
f(a,b,c)=\sqrt{\dfrac{\pi}{a}}e^{\frac{\hat{R}^{2}}{4a}}=\sqrt{\dfrac{\pi}{a}}\sum_{n=0}^{\infty}\dfrac{1}{n!}\left( \dfrac{1}{4a}\right)^{n} R_{2n}(b,c)
\end{equation}
provided that

\begin{equation}\begin{split}
& e^{\sqrt{x^{2}+bx+c}}=e^{\hat{R}x}=\sum_{n=0}^{\infty}\dfrac{x^{n}}{n!}\hat{R}^{n},\\
& \hat{R}^{n}=R_{n}(b,c)
\end{split}\end{equation}
where  $R_{n}(b,c)$ are polynomials of the parameter $b, c$.\\

Even though such a polynomials expansion can be obtained using different procedure, we have tested the validity of our ansatz using the following integral definition

\begin{equation}
R_{m}(b,c)=\dfrac{m!}{2\pi}\int_{0}^{2\pi}e^{-im\phi}e^{\sqrt{e^{2i\phi}+be^{i\phi}+c}}d\phi
\end{equation}
which has been used to benchmark the identity \eqref{B}, with the full numerical integration of \eqref{A}. Further comments will be provided elsewhere.\\\\

\end{document}